\documentclass[11pt,reqno]{amsart}
\usepackage{amscd}

\newcommand{\f}{\operatorname{fl}}
\renewcommand{\r}{\operatorname{rd}}
\newcommand{\p}{\partial}
\newcommand{\Z}{\mathbb Z}

\newcommand{\D}{\mathcal D}

\newcommand{\C}{\mathbb C}

\newcommand{\R}{\mathbb R}

\newcommand{\ep}{\varepsilon}
\renewcommand{\phi}{\varphi}

\newcommand{\cs}{\operatorname{cs}}

\newcommand{\ind}{\operatorname{ind}}

\newcommand{\spin}{\,\operatorname{spin}}

\newcommand{\sign}{\operatorname{sign}}

\newcommand{\bcs}{\beta}

\newtheorem{theorem}{Theorem}[section]
\newtheorem{theorem1}{Theorem}
\newtheorem{lemma}[theorem]{Lemma}

\newtheorem{proposition}[theorem]{Proposition}
\newtheorem{corollary}[theorem]{Corollary}%
\theoremstyle{definition}
\newtheorem{remark}[theorem]{Remark}
\newtheorem{example}[theorem]{Example}
\title{Dirac operators on manifolds with \\ periodic ends}
\thanks{The first author was partially supported by NSF Grant
0505605. The second author was partially supported by NSF
Grant 0305946.}
\author[Daniel Ruberman]{Daniel Ruberman}
\address{Department of Mathematics, MS 050\newline\indent Brandeis
University \newline\indent Waltham, MA 02454}
\email{\rm{ruberman@brandeis.edu}}
\author[Nikolai Saveliev]{Nikolai Saveliev}
\address{Department of Mathematics\newline\indent
University of Miami \newline\indent PO Box 249085
\newline\indent Coral Gables, FL 33124}
\email{\rm{saveliev@math.miami.edu}}
\date\today

\begin{document}
\begin{abstract}
This paper studies Dirac operators on end-periodic spin manifolds
of dimension at least $4$. We provide a necessary and sufficient
condition for such an operator to be Fredholm for a generic
end-periodic metric. We make use of end-periodic Dirac operators
to give an analytical interpretation of an invariant of
non-orientable smooth $4$-manifolds due to Cappell and Shaneson.
 From this interpretation we show that some exotic non-orientable
$4$-manifolds do not admit a metric of positive scalar curvature.
\end{abstract}
\maketitle

\section{Introduction}
Let $M$ be a connected smooth spin $n$--manifold with $n\ge 4$.
With any choice of Riemannian metric $g$ on $M$ one associates the
Dirac operator $D(M,g): C^{\infty} (M,S_g) \to C^{\infty} (M,S_g)$
on the sections of the spinor bundle $S_g$. This is a self--adjoint
elliptic differential operator of order one. Deciding whether this
operator is Fredholm is an important but difficult problem if $M$
is not compact.

In this paper, we study this problem for non--compact manifolds $M$
whose ends are periodic in the sense of Taubes~\cite{T}. In the
first part of the paper, we prove that certain conditions on the
topology of $M$ are sufficient for the Dirac operator $D(M,g)$ to
be Fredholm for a generic periodic metric $g$. The second part
(Sections \ref{S:lift} and \ref{S:psc}) is concerned with
applications of such Dirac operators to $4$-dimensional topology.
We define a metric-dependent invariant of a homology $S^1 \times
S^3$ which lifts its Rohlin
invariant~\cite{ruberman-saveliev:survey}. This leads to a new
analytic
interpretation of an invariant defined by Cappell and Shaneson for
some non-orientable $4$-manifolds in \cite{cappell-shaneson:rp4}
(a different analytic interpretation can be found in
\cite{stolz:spectral}). We then use this interpretation to show
that some of the non-orientable $4$-manifolds constructed in
\cite{cappell-shaneson:rp4} and \cite{akbulut:fake} do not admit
a Riemannian metric of positive scalar curvature. It appears that
this phenomenon is special to dimension four.

A different approach to analyzing the Dirac operator on
non--compact manifolds $M$ is due to Gromov and Lawson \cite{GL}
who showed that $D(M,g)$ is Fredholm if the metric $g$ has positive
scalar curvature at the end of $M$. This analysis is sufficient
for some applications discussed in the second part of our paper;
unfortunately, it does not apply in many situations of geometric
interest. In fact, we were lead to study Dirac operators on
$4$-manifolds with periodic ends by our earlier studies
\cite{ruberman-saveliev:mappingtori,ruberman-saveliev:survey,
ruberman-saveliev:4tori} of gauge theory on manifolds with first
cohomology isomorphic to $\Z$; this paper provides sufficient
conditions, valid in all dimensions, for these operators to be
Fredholm.

To state the result more precisely, recall that any smooth closed
connected spin $n$--manifold $X$ endowed with a smooth map $f: X
\to S^1$ defines an element $[X]\in \Omega^{\spin}_n (S^1)$. Choose
a Riemannian metric $g$ on $X$ and let $\alpha(X)\in KO_n (S^1)$ be
the image of $[X]$ under the homomorphism $\alpha: \Omega^{\spin}_n
(S^1) \to KO_n (S^1)$ described in detail in Section \ref{S:alpha}.
Suppose that $f: X \to S^1$ induces an epimorphism $f_*: \pi_1 X
\to \mathbb Z$ and let $\tilde X$ be the infinite cyclic cover
of $X$ determined by $f_*$ with the induced metric and spin
structure. If $M$ is a spin $n$--manifold with periodic end
modeled on $\tilde X$, any metric on $M$ which agrees with the
metric induced by $g$ over the end will again be called $g$.

\begin{theorem1}\label{T1}
Let $M$ be a spin $n$--manifold with periodic end modeled on $\tilde
X$. Then $\alpha(X) = 0$ if and only if, for a generic metric $g$
on $X$, the Dirac operator $D(M,g): L^2_1 (M,S_g) \to L^2 (M,S_g)$
is Fredholm.
\end{theorem1}

The precise meaning of word ``generic'' is explained in Remark
\ref{R:maier}. Theorem \ref{T1} is deduced from Taubes' fundamental
work on end--periodic differential operators~\cite{T} and the
following result. For any $c \in \mathbb R$, define the twisted
Dirac operator
\[
D^c = D + ic\,f^*(d\theta): C^{\infty}(X,S_g)\to C^{\infty}(X,S_g)
\]
where $d\theta$ is the volume form on $S^1$, and its pull back
$f^*(d\theta)$ acts on the sections of $S_g$ via Clifford
multiplication.

\begin{theorem1}\label{T2}
Let $X$ be a smooth closed connected spin $n$--manifold, $n \ge 4$,
and $f: X \to S^1$ a smooth map such that $f_*: \pi_1 X \to \Z$ is
onto. Then $\alpha (X)\in KO_n (S^1)$ vanishes if and only if, for
a generic metric $g$ on $X$, we have $\ker D^c (X,g) = 0$ for all
$c \in \mathbb R$.
\end{theorem1}

\noindent\textbf{Acknowledgments:} We thank Lev
Kapitanski, Tom Mrowka, Jonathan Rosenberg, and Peter Teichner
for helpful conversations on the ideas discussed in this paper.  We
also thank Ulrich Bunke for communicating to us a proof (see
Section~\ref{S:vanish}) of the `if' direction in Theorem~\ref{T2}.  Parts of this project were accomplished while we attended the workshop on $4$-manifolds at Oberwolfach and (the second author) the Summer Session of the Park City/IAS Mathematics Institute. We express our appreciation to
the organizers of both for providing a stimulating environment.


\section{Deducing Theorem \ref{T1} from Theorem \ref{T2}}
\label{S:taubes}
Let $X$ be a connected closed smooth spin $n$--manifold with $n \ge
4$. Choose a Riemannian metric $g$ on $X$ and associate with it the
self--adjoint Dirac operator $D(X,g): C^{\infty}(X,S_g) \to
C^{\infty} (X,S_g)$ on the sections of the spinor bundle $S_g$.

Given a smooth map $f: X \to S^1$ inducing an epimorphism $f_*:
\pi_1 X \to \Z$,  consider the infinite cyclic covering $\tilde X
\to X$ associated with $f_*$. Then $\tilde X$ inherits a Riemannian
metric from $X$, called again $g$, and a spin structure, which in
turn give rise to the periodic Dirac operator $D (\tilde X,g):
C^{\infty}(\tilde X,S_g)\to C^{\infty}(\tilde X,S_g)$. We wish to
prove that, for a generic metric $g$ on $X$, the $L^2$--completion
of this operator, $D(\tilde X,g): L^2_1\,(\tilde X,S_g)\to L^2\,
(\tilde X,S_g)$, is Fredholm. The Fredholmness of more general Dirac
operators claimed in Theorem \ref{T1} will then follow by the
excision principle, see Lockhart--McOwen \cite{LM}. The following
construction is due to Taubes \cite{T}.

Observe that $S_g \to \tilde X$ is a pull back to $\tilde X$ of the
spinor bundle $S_g \to X$. Let $\tau: \tilde X \to \tilde X$ be a
deck transformation and use the same symbol $\tau$ to denote its
lift to an automorphism of the bundle $S_g \to \tilde X$. Choose a
smooth map $\tilde f: \tilde X\to \mathbb R$ lifting $f: X\to S^1$.

Given $\psi \in C^{\infty}_0 (\tilde X,S_g)$ and $z \in \mathbb C^*$,
denote by $\tau^*(\psi)$ the pull back of $\psi$ and define
\[
\hat\psi_z = z^{\tilde f}\cdot\sum_n z^n \cdot (\tau^*)^n (\psi),
\]
for a fixed choice of branch of $\ln z$. Then $\hat\psi_z$ is a smooth
section of $S_g \to \tilde X$ such that
\[
\tau^*(\hat\psi_z) = z^{\tilde f + 1}\cdot\sum_n z^n\cdot(\tau^*)^{n+1}(\psi)
= z^{\tilde f}\cdot\sum_n z^{n+1} (\tau^*)^{n+1}(\psi) = \hat\psi_z,
\]
therefore, $\hat\psi_z$ defines a smooth section of $S_g \to X$. This
operation is referred to as the Fourier--Laplace transform. It converts
the Dirac operator $D(\tilde X,g)$ into the family of operators
\[
D_z (X,g) = z^{\tilde f}\cdot D(X,g)\cdot z^{-\tilde f} = D(X,g) +
z^{\tilde f}\cdot [D(X,g),z^{-\tilde f}],
\]
where $D(X,g): C^{\infty} (X,S_g) \to C^{\infty} (X,S_g)$ is the Dirac
operator on $X$, and $z^{\tilde f} \cdot [D(X,g),z^{-\tilde f}]$ is a
zero order operator (a bundle automorphism of $S_g$). The following is
a special case of the more general result of Taubes \cite[Lemma 4.3]{T}.

\begin{theorem}\label{T:taubes}
The operator $D(\tilde X,g): L^2_1\,(\tilde X,S_g)\to L^2\,(\tilde X,
S_g)$ is Fredholm if and only if $\ker D_z (X,g) = 0$ for all $z \in
\C^*$ with $|z|= 1$.
\end{theorem}

A direct calculation shows that, if $z = e^{-ic}$ for $c\in \R$, then
$D_z (X,g) = D (X,g) + ic\,f^*(d\theta)$, where $d\theta$ is the
standard volume form on $S^1$. These are precisely the twisted Dirac
operators $D^c (X,g)$ of Theorem \ref{T2}.


\section{Transporting invertibility}\label{S:transport}
Let $X$ and $X'$ be connected closed smooth spin $n$--manifolds with
$n\ge 4$. We will say that $f: X \to S^1$ is \emph{spin cobordant
to $f': X'\to S^1$ over $S^1$} if they define the same element in
the group $\Omega^{\spin}_n (S^1)$. This means that one can find an
oriented smooth spin cobordism $W$ from $X$ to $X'$ and a smooth
function $F: W\to S^1$ restricting to $f$ and $f'$ respectively at
the two boundary components. The main technical step in proving
Theorem \ref{T2} is the following result.

\begin{theorem}\label{T:main}
Suppose that $f: X \to S^1$ is spin cobordant to $f': X' \to S^1$
over $S^1$ via a cobordism admitting a handle decomposition (starting
at $X$) with
handles of index at most $n - 1$. If there is a metric $g$ on $X$
such that $\ker D^c (X,g) = 0$ for all $c$ then one can find a
metric $g'$ on $X'$ such that $\ker D^c (X',g') = 0$ for all $c$.
\end{theorem}

\begin{remark}\label{R:maier}
Finding just one metric $g$ as in the above theorem is sufficient
to show that such metrics form a non--empty subset of the space of
Riemannian metrics on $X$ which is open in $C^1$--topology and
dense in all $C^k$--topologies, $k \ge 1$, see Maier
\cite[Proposition 3.1]{maier}. The metrics in the above subset are
referred to as ``generic''.
\end{remark}

The proof of Theorem \ref{T:main} closely follows the proof of
Theorem 1.2 in the paper of Ammann, Dahl, and Humbert~\cite{ADH},
which is a stronger assertion for the (untwisted) Dirac
operator $D(X,g)$: they prove existence of metric $g'$ on $X'$
such that $\dim\ker D(X',g')\le \dim\ker D(X,g)$ even when $\ker
D(X,g) \ne 0$. The analytical estimates that are at the heart of
the proof of Theorem \ref{T:main} are very similar to those which
lie behind~\cite[Theorem 1.2]{ADH}, except that one has to pay
extra attention to the additional parameter $c$. For this reason,
this section will provide a sketch that consists of pointers to
the relevant lemmas in~\cite{ADH} together with comments on the
modifications necessary to take account of $c$. We will be happy
to supply full details to the interested reader.

The proof of Theorem \ref{T:main} will proceed by induction on the
number of handles in the spin cobordism.  Adding a 0--handle amounts
to taking the disjoint union with an $(n+1)$--ball, which does not
raise $\dim\ker D^c(X,g)$ since the $n$-sphere has positive scalar
curvature. So we may assume in our induction that we are adding
handles of index at least $1$.

Let $W$ be a spin cobordism over $S^1$ from
$f: X \to S^1$ to $f': X' \to S^1$ with just one handle of index
$k + 1$, $0 \le k \le n - 2$, so that
\[
W = ([0,1]\times M) \cup (D^{\,k+1} \times D^{\,n-k}),
\]
where the handle is glued along $S^k \times D^{n-k} \subset X$.
Then $X'$ can be obtained from $X$ by cutting out $S^k \times
D^{n-k}$ and gluing back $D^{k+1}\times S^{n-k-1}$\,:
\begin{equation}\label{E:surgery}
X' = (X - S^k\times D^{\,n-k})\cup (D^{\,k+1}\times S^{n-k-1}).
\end{equation}
In this $n$--dimensional description, the fact that $W$ is a spin
cobordism over $S^1$ means that the spin--structure on $X'$ and
the smooth function $f': X'\to S^1$ agree with the spin structure
on $X$ and the function $f: X \to S^1$ over $X - (S^k \times
D^{\,n-k}) = X' - (D^{\,k+1} \times S^{n-k-1})$.

Theorem \ref{T:main} will follow as soon as we show that, if
$X'$ is obtained from $X$ by operation (\ref{E:surgery}) and
$\ker D^c(X,g) = 0$ for all $c$, then there exists a metric $g'$
on $X'$ such that $\ker D^c(X',g') = 0$ for all $c$.


\subsection{Preliminaries}
Let $X$ be a connected closed smooth spin $n$--manifold equipped
with a smooth function $f: X\to S^1$, and let $D^c(X,g) = D(X,g) +
ic\,f^* d\theta$, $c\in\mathbb R$, be the twisted Dirac operators.

\begin{lemma}\label{L:per}
For any real number $c$, the spectra of $D^{c+1}(X,g)$ and $D^c(X,g)$
coincide. In particular, $\ker D^{c + 1}(X,g) = \ker D^c(X,g)$.
\end{lemma}

\begin{proof}
View $\theta$ as a multivalued function on $S^1$ and observe that
$e^{i f^*\theta}$ is a smooth function on $X$. Then
$e^{-if^*\theta} D^c(X,g)\,(e^{i f^*\theta}\psi) =
e^{-if^*\theta} (de^{if^*\theta}\psi + e^{if^*\theta} D^c(X,g)\psi)
= i f^* d\theta \psi + D^c(X,g)\psi = D^{c+1}(X,g)\psi$ so that
$D^{c + 1}(X,g)$ and $D^c(X,g)$ have the same spectrum.
\end{proof}

\begin{lemma}\label{L:spec}
If the form $f^* d\theta$ is exact then, for any real $c$, the
spectrum of $D^c(X,g)$ equals the spectrum of $D(X,g)$.
\end{lemma}

\begin{proof}
Write $f^* d\theta = du$ for some smooth function $u: X \to
\mathbb R$. Then $D^c(X,g) = D(X,g) + ic\,du = e^{-icu}\, D(X,g)\;
e^{icu}$, therefore, the spectra of $D^c(X,g)$ and $D(X,g)$ are
the same for all $c\in\mathbb R$.
\end{proof}

In other language, these lemmas say that the operators $D^c$ and
$D^{c+1}$ (resp. $D$ and $D^c$ when $f^*d\theta$ is exact) are gauge
equivalent.

\begin{corollary}\label{C:est}
Let $S^{\ell}$, $\ell\ge 1$, be an $\ell$--dimensional sphere
with round metric $g^{\r}$ of radius one and bounding spin
structure, and $M = N \times S^{\ell}$ a closed spin manifold
with a product spin structure and a product metric $g = g_N +
g^{\r}$. For any $f: M \to S^1$ such that $f^*d\theta$ is exact,
the spectrum of $D^c(M,g)^2$ is bounded from below by $\ell^2/4$.
\end{corollary}

\begin{proof}
Let $f^* d\theta = du$ then the proof of Lemma \ref{L:spec} shows
that $D^c(M,g) = e^{-icu}\, D(M,g)\; e^{icu}$ therefore
$D^c(M,g)^2 = e^{-icu}\,D(M,g)^2\;e^{icu}$ and $D^c(M,g)^2$ and
$D(M,g)^2$ have the same spectrum. The rest of the proof is as in
\cite[Proposition 2.5]{ADH}, using estimates on the spectrum of
$D(S^{\ell},g^{\r})^2$. If $\ell \ge 2$, this spectrum is bounded
from below by $\ell^2/4$, see for instance \cite{B}. If
$\ell = 1$ and the spin structure extends over $D^2$, we have
$D(S^1,g^{\r}) = i\,d/d\theta + 1/2$, which obviously gives the
desired estimate.  (Note that, with respect to the non--bounding
spin structure, the Dirac operator is of the form $D(S^1,g^{\r})
= i\,d/d\theta$ and hence has a non--zero kernel.)
\end{proof}

In the analysis of the Dirac operator on a spin manifold, a key
role is played by the Lichnerowicz formula for its square.  For
a twisted Dirac operator, this formula has an extra term involving
the curvature of the twisting connection, see for instance
\cite[page 134]{BGV}. In our case the twisting connection $ic
f^* d\theta$ is flat, so the formula becomes
\begin{equation}\label{E:lich}
D^c(X,g)^2 = (\bar\nabla^{g,c})^*\,\bar\nabla^{g,c}\; +
\; \frac {\kappa} 4,
\end{equation}
where $\kappa$ is the scalar curvature of $(X,g)$ and
$\bar\nabla^{g,c}$ is the induced (twisted) connection on the
spin bundle $S_g$.

\medskip

Our final observation is that, for any metrics $g$ and $g'$ on $X$,
there exists a unique automorphism $b: TX\to TX$ which is positive,
symmetric with respect to $g$, and has the property that $g
(\xi,\eta) = g'(b(\xi),b(\eta))$. The map on orthonormal frames
induced by $b$ gives rise to a map $\beta^g_{g'}: S_g \to S_{g'}$
of spinor bundles associated with the metrics $g$ and $g'$. Note
that this map preserves the fiberwise length of spinors and
that $\beta^{g'}_g\circ\beta^g_{g'} = 1$. The diagram
\[
\begin{CD}
C^{\infty}(X,S_g) @> \beta^g_{g'} >> C^{\infty}(X,S_{g'}) \\
@V D^c (X,g) VV  @VV D^c (X,g') V \\
C^{\infty}(X,S_g) @> \beta^g_{g'} >> C^{\infty}(X,S_{g'}).
\end{CD}
\]

\bigskip\noindent
need not commute; however, there is an explicit formula relating
the operators $D^c(X,g)$ and $(\beta^g_{g'})^{-1}\circ D^c(X,g')
\circ\beta^g_{g'}: C^{\infty}(X,S_g)\to C^{\infty}(X,S_g)$.  All
we need is its rather general form; compare with \cite{BG} or
\cite{maier},
\begin{equation}\label{E:bg}
((\beta^g_{g'})^{-1}\circ D^c(X,g')\circ\beta^g_{g'})(\psi) =
D^c(X,g)(\psi) + A(\bar\nabla^{g,c}\psi) + B(\psi),
\end{equation}
where $A: T^*X \otimes S_g \to S_g$ and $B: S_g \to S_g$ are
bundle maps independent of $c$ such that
\begin{equation}\label{E:AB}
|A| \le C\,|g-g'|_g\quad\text{and}\quad
|B| \le C\,(|g-g'|_g + |\nabla^g (g-g')|_g)
\end{equation}
for some constant $C$. In the above formulas, $\nabla^g$ is the
Levi--Civita connection on $TX$ associated with metric $g$, and
$\bar\nabla^{g,c}$ is the induced (twisted) connection on $S_g$.


\subsection{Approximating by a product metric}
Let $g$ be a metric on $X$ and denote by $h$ the induced metric on
$S^k \subset X$, $0 \le k \le n - 2$. Use the exponential map to
identify a tubular neighborhood $U(R)$ of $S^k$ of radius $R > 0$
with the product $S^k \times D^{n-k}$. This neighborhood has two
metrics\,: one is the original metric $g$ and the other the product
metric $h + g^{\f}$, where $g^{\f}$ is the flat metric on $D^{n-k}$.

Choose a small real number $\delta > 0$ and a smooth cut--off
function $\eta: X \to \mathbb R$ such that $0 \le \eta \le 1$,
$\eta = 1$ on $U(\delta)\subset X$, $\eta = 0$ on $X - U(2\delta)$,
and $|d\eta|_g \le 2/\delta$. Consider the metric
\[
g_{\delta} = \eta\,(h + g^{\f}) + (1 - \eta)\,g,
\]
which is a product metric in $U(\delta)$. The proof of the following
proposition, which is along the lines of~\cite[Proposition 3.3]{ADH}, will be
sketched in the remainder of this subsection.

\begin{proposition}\label{P:prod}
Let $g$ be a metric on $X$ such that $\ker D^c(X,g) = 0$ for all $c$.
Then, if $\delta > 0$ is sufficiently small, $\ker D^c(X,g_{\delta})
= 0$ for all $c$.
\end{proposition}

Suppose this is not true. Then one can find a sequence $\delta_m \to
0$ and a sequence $c_m$ such that $\ker D^{c_m}(X,g_{\delta_m}) \ne
0$. Lemma \ref{L:per} implies that there is no loss of generality
in assuming that $0\le c_m\le 1$ for all $m$. Therefore, passing to
a subsequence if necessary, we will assume that $c_m$ is a convergent
sequence and denote its limit by $c$.

For the sake of simplicity, we will use the notations $g_m =
g_{\delta_m}$ and $\beta^{g_m}_g = \beta_m$. Choose a non--zero
spinor $\phi_m \in \ker D^{c_m}(X,g_m)$ and normalize it so that
\begin{equation}\label{E:norm1}
\int_X |\phi_m|^2\,dv^{g_m} = 1.
\end{equation}

\begin{lemma}
The sequence $\beta_m\phi_m$ is bounded in $L^2_1\,(X,g)$.
\end{lemma}

\begin{proof}
Observe that the estimates (\ref{E:bg}) are uniform in $c$, hence the
proof of \cite[Lemma 3.4]{ADH} generalizes to our case word for word
to show that the sequences $\beta_m\phi_m$ and $\bar\nabla^{g,c_m}
(\beta_m\phi_m)$ are bounded in $L^2(X,g)$. To finish the proof,
observe that
\[
\bar\nabla^{g,c_m}(\beta_m\phi_m) = \bar\nabla^g(\beta_m\phi_m) +
ic_m f^*d\theta\,\otimes\,\beta_m \phi_m
\]
hence
\[
|\bar\nabla^g(\beta_m\phi_m)|_g^2\quad\le\quad 2\,|\bar\nabla^{g,c_m}
(\beta_m\phi_m)|_g^2\; + \;2\,c_m^2\,|f^*d\theta|_g^2\,|\beta_m\phi_m|^2
\]
with $c_m^2 \le 1$. Therefore,
\[
\int_X\;|\bar\nabla^g(\beta_m\phi_m)|_g^2\,dv^g\; \le\; C\,\int_X
\;|\bar\nabla^{g,c_m}(\beta_m\phi_m)|_g^2\,dv^g\; +\;
C\int_X\;|\beta_m\phi_m|^2\,dv^g,
\]
and the sequence
\[
\int_X\;|\bar\nabla^g(\beta_m\phi_m)|_g^2\,dv^g
\]
is bounded.
\end{proof}

Since $\beta_m\phi_m$ is bounded in $L^2_1 (X,g)$, after passing to a
subsequence if necessary, we may assume that $\beta_m\phi_m$ converges
to a spinor $\phi$ weakly in $L^2_1(X,g)$ and strongly in $L^2(X,g)$.
Because of the strong convergence in $L^2(X,g)$, the normalization
(\ref{E:norm1}) is preserved in the limit; in particular, $\phi\ne 0$.
To continue, we will need the following extension of \cite[Lemma 2.2]
{ADH}.

\begin{lemma}\label{L:cpt}
Let $(M,g)$ be a Riemannian manifold, not necessarily compact, and $f:
M \to S^1$ a continuous map. Suppose that $D^c(M,g)(\phi) = 0$ for
some $c \in [0,1]$. Then, for any compact set $K \subset M$, there
exists a constant $C = C(K,M,g)$ independent of $c\in [0,1]$ such that
\[
\|\phi\|_{\,C^2(K,g)}\;\le\; C\,\|\phi\|_{\,L^2(M,g)}.
\]
\end{lemma}

\begin{proof}
This follows by applying a standard bootstrapping argument to the equation
$D^c(M,g)(\phi) = 0$. More precisely, write this equation in the form
$D(M,g)(\phi) = -ic f^*(d\theta)\,\phi$ and apply $D(M,g)$ to both
sides of it to obtain
\[
D(M,g)^2(\phi) = c^2\, |f^*(d\theta)|^2_g\;\phi.
\]
Since $0 \le c \le 1$, the $L^2(M,g)$--norm of the right hand side can
be estimated from above by $C\,\|\phi\|_{L^2(M,g)}$, where $C$ does not
depend on $c$. Therefore, for any compact set $K_1\subset M$, there is
a constant $C$ independent of $c$ such that
\[
\|\phi\|_{\,L^2_2(K_1,g)}\;\le\;C\,\|\phi\|_{\,L^2(M,g)}.
\]
Repeat the argument to obtain the estimate
\[
\|\phi\|_{\,L^2_4(K_1,g)}\;\le\;C\,\|\phi\|_{\,L^2(M,g)}.
\]
If the boundary of $K_1$ is smooth, then we have the Sobolev embedding
$L^2_4(K_1,g) \to C^1 (K_1,g)$ so that $\|\phi\|_{C^1(K_1,g)}\;\le\; C
\,\|\phi\|_{L^2(X,g)}$ with a constant $C$ independent of $c$. Now one
can use Schauder's estimates as in \cite[Lemma 2.2]{ADH} to obtain the
$C^2$--estimate.
\end{proof}

Apply the above lemma to the sequence $\beta_m\phi_m$ to conclude that,
for any sufficiently small $\ep > 0$, the sequence $\beta_m\phi_m$ is
bounded in $C^2(X - U(\ep))$. The Ascoli Theorem then implies that
a subsequence of $\beta_m\phi_m$ converges in $C^1(X - U(\ep))$. The
limit $\phi$ is in $C^1_{\operatorname{loc}}(X - S)$ and satisfies the
equation $D^c(X,g)(\phi) = 0$ on $X - S$. The removal of singularities
lemma \cite[Lemma 2.4]{ADH} says that $\phi$ extends to a smooth
spinor on $X$ such that $D^c(X,g)(\phi) = 0$. This contradicts the
assumption that $\ker D^c(X,g) = 0$ for all $c$, and proves
Proposition~\ref{P:prod}.


\subsection{Metrics which are product in the surgery region}
Let $X'$ be obtained from $X$ by surgery along an embedded sphere
$S^k$ as in (\ref{E:surgery}). According to Proposition~\ref{P:prod},
we may assume that the metric on $X$ is a product on a tubular
neighborhood $U(R)$ of $S^k$ of radius $R > 0$. For any sufficiently
small number $\rho > 0$, Ammann, Dahl and Humbert \cite[Section 3.2]
{ADH} construct a metric $g_{\rho}$ on $X'$ which coincides with $g$
on $X - U(R)$, is conformally equivalent to $g$ on the entire region
$X - U(2\rho)$, and is a product metric $\gamma_{\rho} + g^{\r}$ on
the region $D^{k+1} \times S^{n-k-1}$ glued into $X - U(\rho/2)$ by
the surgery.

\begin{proposition}\label{P:two}
Suppose that $\ker D^c(X,g) = 0$ for all $c$. Then, for any metric
$g_{\rho}$ with sufficiently small $\rho > 0$, we have $\ker D^c
(X',g_{\rho}) = 0$ for all $c$.
\end{proposition}

Our proof follows rather closely the proof in~\cite[Section 3.2]{ADH},
with some additional observations. Suppose that the proposition is
not true. Then one can find a sequence $\rho_m > 0$ converging to
zero and a sequence $c_m \in \mathbb R$ such that $\ker D^{c_m}
(X',g_{\rho_m})\ne 0$ for all $m$.

Lemma \ref{L:per} implies that there is no loss of generality in
assuming that $0\le c_m\le 1$ for all $m$. Therefore, passing to a
subsequence if necessary, we will assume that $c_m$ is a convergent
sequence and denote its limit by $c$.

For each $m$, choose a spinor $\psi_m \in \ker D^{c_m}(X',g_{\rho_m})$
and view it as a spinor on $X' - U'(2\rho_m) = X - U(2\rho_m)$, two
manifolds with conformally equivalent metrics. Since the kernel of
a (twisted) Dirac operator is conformally invariant, $\psi_m$ gives
rise to a spinor $\phi_m$ in the kernel of $D^{c_m}(X,g)$ on $X -
U(2\rho_m)$. Fix a real number $s_m$ slightly greater than $2\rho_m$
and normalize $\phi_m$ so that

\begin{equation}\label{E:norm}
\int_{X - U(s_m)}|\phi_m|^2\,dv^g = \int_{X' - U'(s_m)}|\phi_m|^2
\,dv^g = 1.
\end{equation}

\medskip

\begin{lemma}
The sequence $\phi_m$ is bounded in $L^2 (X - U(\delta))$ for any
choice of $\delta \in (0,R)$ by a constant independent of both
$\delta$ and $c_m$.
\end{lemma}

\begin{proof}
In the untwisted case, this is derived in \cite{ADH} as follows. The
metric $g_{\rho}$ is a product metric on $U'(2s_m) = D^{k+1} \times
S^{n-k-1}$ of the form $h_{\rho} + g^{\r}$. Embed the disk $(D^{k+1},
h_{\rho})$ isometrically into a closed spin Riemannian manifold
$(N,g_N)$ of dimension $k + 1$, and let $M = N \times S^{n-k-1}$
with product metric $g_M = g_N + g^{\r}$ and the bounding spin
structure on $S^{n-k-1}$. The spectrum of $D (M,g_M)^2$ is bounded
from below by $(n-k-1)^2/4$, which is positive because $k \le n - 2$.
This leads to the desired estimate on the $L^2$ norm of $\phi_m$
through the use of the Rayleigh quotients for the operator $D
(M,g_M)^2$ over $U'(2s_m) \subset M$, see Proposition 3.6 and
estimate (17) in \cite{ADH}.

The same approach works for the operators $D^{c_m}(M,g_M)^2$ to give
an estimate which is uniform in $c_m$ as long as we can apply
Corollary \ref{C:est}, or in other words, as long as we know that
$(f')^*(d\theta)\in \Omega^1 (U'(2s_m))$ is the restriction of an
exact 1--form on $M$. The latter can be seen as follows. The
restriction of the map $f': X' \to S^1$ to $U'(2s_m)$ can be
included into the following commutative diagram

\[
\begin{CD}
S^{n-k-1} @>i'>> U'(2s_m) \\
@Vi VV  @VV f' V \\
D^{n-k} @> f >> S^1
\end{CD}
\]

\bigskip\noindent
where the maps $i$ and $i'$ are induced by the inclusions of
$S^{n-k-1} = \p D^{n-k} \subset S^k\times D^{n-k} = U(\rho_m/2)$
and $S^{n-k-1}\subset D^{k+1}\times S^{n-k-1} = U'(\rho_m/2)$,
respectively. The map $f\circ i: S^{n-k-1} \to S^1$ factors
through the contractible $D^{n-k}$ and hence is homotopic to
zero. Together with the fact that $(i')^*: H^1 (U'(2s_m)) \to
H^1 (S^{n-k-1})$ is an isomorphism, this implies that $(f')^*:
H^1(S^1) \to H^1(U'(2s_m))$ is zero. In particular,
$(f')^*(d\theta)$ on $U'(2s_m)$ is exact (of course, this form
is automatically exact unless $k = n - 2$). We let $(f')^*(d\theta)
= du'$ for some $u': U'(2s_m) \to \mathbb R$ and extend $u'$ to a
smooth function $u$ on $M = N \times S^{n-k-1}$.
\end{proof}

Let $N$ be a sufficiently large positive integer. Then the sequence
$\phi_m$ is bounded in $L^2(X - U(1/N))$ and, by Lemma \ref{L:cpt},
also in $C^2(X - U(2/N))$. By Ascoli's Theorem, there is a
subsequence $\phi^0_m$ which converges in $C^1(X - U(2/N))$ to a
spinor $\phi^0$. By induction construct, for each $i\ge
1$, a subsequence $\phi^i_m$ of $\phi^{i-1}_m$ converging to
a spinor $\phi^i$ in $C^1(X - U(2/(N + i)))$. Obviously,
$\phi^i$ is an extension of $\phi^{i-1}$ for $i \ge 1$. Use the
diagonal algorithm to construct a sequence $\phi'_m$ converging
in $C^1_{\operatorname{loc}}(X - S^k)$ to a spinor $\phi'$.

As $D^{\,c_m}(X,g)(\phi'_m) = 0$ over $X - U(2\rho_m)$ with $\lim
\rho_m = 0$ and $\lim c_m = c$, the $C^1_{\operatorname{loc}}
(X - S^k)$ convergence implies that $D^c(X,g)(\phi') = 0$ on $X
- S^k$. The removal of singularities lemma \cite[Lemma 2.4]{ADH}
says that $\phi'$ extends to a smooth spinor on all of $X$ such
that $D^c(X,g)(\phi') = 0$. Moreover, the normalization
(\ref{E:norm}) is preserved in the limit thus implying that $\phi'
\ne 0$. This contradicts the assumption that $\ker D^c (X,g) = 0$
for all $c$.  Proposition~\ref{P:two} follows.


\section{Proof of Theorem \ref{T2}}\label{S:alpha}
In this section, we first give a detailed description of the
invariant $\alpha$, and then use the results of Section
\ref{S:transport} on transporting the invertibility of the
twisted Dirac operators across a cobordism to prove Theorem \ref{T2}.


\subsection{The invariant $\alpha$}\label{S:alpha1}
Every smooth closed spin $n$--manifold $X$ equipped with a map
$f: X \to B\pi$ to the classifying space of a discrete group
$\pi$ defines an element $[X]$ in the spin cobordism group
$\Omega^{\spin}_n (B\pi)$. By definition, $\alpha (X) \in KO_n
(B\pi)$ is the image of $[X]$ under the natural transformation
of homology theories $\alpha: \Omega^{\spin}_* \to KO_*$
defined in Milnor \cite{milnor}.

The only two cases we will be interested in are those of $\pi = \Z$
(corresponding to spin manifolds over a circle) and $\pi = \{1\}$
(corresponding to spin manifolds). As a matter of convenience,
the respective $\alpha$--invariants will be denoted by $\alpha(X)
\in KO_n (S^1)$ and $\alpha_n (X) \in KO_n$. They can be described
explicitly as follows, compare with B\"ar--Dahl \cite[Section 3]
{bar-dahl}.

A straightforward calculation shows that $KO_n (S^1) = KO_n \,
\oplus\,KO_{n-1}$. Then, for any spin $n$--manifold $X$ endowed
with a smooth map $f: X \to S^1$,
\[
\alpha (X) = \alpha_n (X) + \alpha_{n-1} (Y) \in KO_n
\,\oplus\,KO_{n-1},
\]
where $Y = f^{-1} (p)$ for any choice of a regular value $p\in S^1$.
Moreover, $KO_k = \Z$ if $k = 0, 4 \pmod 8$, $KO_k = \Z/2$ if $k =
1, 2 \pmod 8$, $KO_k = 0$ otherwise, and
\[
\alpha_k (X) =
\begin{cases}
\;\ind D^+ (X,g) & \text{if $n = 0 \pmod 8$,} \\
\;1/2\,\ind D^+ (X,g) & \text{if $n = 4 \pmod 8$,} \\
\;\dim \ker D (X,g)\pmod 2 & \text{if $n = 1 \pmod 8$,} \\
\;\dim \ker D^+ (X,g)\pmod 2 & \text{if $n = 2 \pmod 8$.}
\end{cases}
\]

\medskip\noindent
where $D^+ (X,g): C^{\infty}(X,S^+_g) \to C^{\infty}(X,S^-_g)$
is the chiral Dirac operator. In particular, $\alpha_4 (X) =
- \sign (X)/16$ for any spin $4$-manifold $X$.

It is clear from the above that the invertibility of $D (X,g)$
on a spin $n$--manifold $X$ implies that $\alpha_n (X) = 0 \in
KO_n$.

\subsection{Vanishing of $\alpha(X)$ is necessary}\label{S:vanish}
We wish to prove that, if $X$ is a spin manifold of dimension
$n \ge 4$ over a circle which admits a metric $g$ such that
the operators $D^c (X,g)$ are invertible for all
$c \in \mathbb R$, then $\alpha (X) = 0 \in KO_n (S^1)$.

This is certainly true for $n = 0, 4, 6, 7\pmod 8$ because, in
these dimensions, $KO_{n-1} = 0$ so that $\alpha (X) = \alpha_n
(X)$, and the latter vanishes for invertible $D (X,g)$ (in
particular, the invertibility of all $D^c (X,g)$ is not really
needed). The proof of the general case sketched below was
suggested to us by Ulrich Bunke.

Let $C^*_{\R} (\Z)$ be the real group $C^*$--algebra of $\Z$. It
follows from the positive solution to the real version of the
Baum--Connes conjecture for group $\Z$ that the assembly map $A:
KO_n (S^1) \to KO_n (C^*_{\R}(\Z))$ is injective. Therefore, it
is sufficient for us to show that $A(\alpha(X))$ vanishes.

The class $A(\alpha(X)) \in KO_n (C^*_{\R}(\Z))$ admits the
following description. The involution $c \to -c$ combined
with the real structure on $D (X,g)$ makes operators $D^c (X,g)$
into a real family $\D$ acting on a Hilbert bundle
over $S^1$. The space of sections of this bundle is naturally a
Hilbert module over $C(S^1)$. The Fourier transform makes it
into a module over $C^*(\Z)$ and further into a real module over
$C^*_{\R}(\Z)$ if one takes into account the above real structure,
as in the paper of Bunke--Schick \cite{bunke-schick}. The $K$--theoretic
index of the family $\D$ is then $A(\alpha(X)) \in KO_n (C^*_{\R}
(\Z))$.

In this interpretation, the invertibility of all $D^c (X,g)$
means that the family $\D$ is invertible on each fiber over
$S^1$. Since $S^1$ is compact, one can find a positive uniform
bound from below on the spectrum of $D^c (X,g)^2$, which makes
the family of inverses uniformly bounded. This implies of
course that the index of $\D$ is zero so that $A(\alpha(X))$
vanishes.


\subsection{Vanishing of $\alpha(X)$ is sufficient in dimension 4}
A direct calculation (using for instance the Atiyah--Hirzebruch
spectral sequence) shows that $\Omega^{\spin}_4 (S^1) =
\Omega^{\spin}_4\,\oplus\,\Omega^{\spin}_3$, and it is well known
that $\Omega^{\spin}_3 = 0$ and $\Omega^{\spin}_4$ is isomorphic
to $\Z$ via $[X] \to \sign(X)/16$. Therefore, if $\alpha (X) =
-\sign (X)/16 = 0$, there is a spin cobordism over a circle from
$S^4$ to $X$. The sphere $S^4$ has a metric $g$ of positive scalar
curvature, hence all operators $D^c (X,g)$ are invertible by
(\ref{E:lich}). Theorem \ref{T2} is now a consequence of Theorem
\ref{T:main} and the following Lemma (which we prove in arbitrary
dimension).

\begin{lemma}\label{L:handles}
If $(X,f)$ of dimension $n\ge 4$ is spin cobordant over a circle
to $(X',f')$ such that $f'_*: \pi_1 X' \to \mathbb Z$ is onto then
$(X,f)$ is spin cobordant over a circle to $(X',f')$ via a cobordism
without $n$--handles.
\end{lemma}

\begin{proof}
Any spin cobordism over a circle from $(X,f)$ to $(X',f')$ can be
split into a union $U \cup V$, where $U$ consists of $k$--handles,
$1 \le k \le n-1$, and $V$ of $n$--handles.  Both $U$ and $V$ are
spin cobordisms over a circle glued along  $(X'',f'')$, where $X''
= X'\,\#\,\ell\,(S^1 \times S^{n-1})$ and $f'': X'' \to S^1$
coincides on the punctured copy of $X'$ with $f': X' \to S^1$.
This can be seen by flipping the cobordism, so that $V$ is obtained
by attaching 1--handles to $[0,1] \times X'$.  With an extra effort,
one can also ensure that the restriction of $f'': X'' \to S^1$ to
each copy of $S^1 \times S^{n-1}$ is homotopic to zero\,:  just slide
one foot of a 1--handle if necessary around $X'$ and use the fact
that $f_*': \pi_1 X' \to \Z$ is onto.

Next,  construct a cobordism $V'$ from $X''$ to $X'$ by attaching
2--handles to $[0,1] \times X''$, one for each copy of $S^1\times
S^{n-1}$, along embedded circles $S^1 \subset S^1 \times S^{n-1}$
generating $H_1 (S^1 \times S^{n-1}; \Z)$ and having bounding spin
structure (to ensure that the resulting cobordism is spin). Since
the restriction of $f''$ to each of the $S^1 \times S^{n-1}$ in
$X''$ is homotopic to zero, the cobordism $V'$ is automatically
over a circle. The union $U \cup V'$ is then the desired cobordism.
\end{proof}


\subsection{Vanishing of $\alpha(X)$ is sufficient in dimensions
$n \ge 5$}
Let $X$ be a spin manifold of dimension $n\ge 5$ and $f: X\to S^1$
be a smooth map such that $f_*: \pi_1 X \to \Z$ is onto. One can
add 2--handles to $[0,1]\times X$ to get a spin cobordism $W$ to a
spin manifold $X'$ over $S^1$ so that $f_*': \pi_1 X'\to \Z$ is an
isomorphism. Since the 2--handles are attached along circles
mapped to zero by $f_*$, the cobordism $W$ is over a circle.  The
cobordism $W$ can be flipped so that $X$ is obtained from $X'$ by
adding $n-1$ handles.

The existence of $W$ has two important implications.  One is that
$\alpha (X') = \alpha (X) = 0$, and the other is that $D^c (X,g)$
are invertible for a generic metric $g$ on $X$ if and only if
$D^c (X',g')$ are invertible for a generic metric $g'$ on $X'$,
see Theorem \ref{T:main}.

Now, according to Joachim and Schick \cite{JS}, the vanishing of
$\alpha(X')$ implies existence of a metric $g'$ of positive scalar
curvature on $X'$ (this is known as the Gromov--Lawson--Rosenberg
conjecture; it has been proved in particular for manifolds whose
fundamental group is free abelian). Therefore, all operators $D^c
(X',g')$ are invertible by (\ref{E:lich}).

\begin{remark}
Note that the above construction does not necessarily result in a
metric of positive scalar curvature on $X$ itself, because $X$ is
obtained from $X'$ using handles of codimension two. While the
invertibility of $D^c$ can be transported across such a handle,
transporting the positive scalar curvature condition would require
codimension at least three.
\end{remark}


\section{An integral lift of the Rohlin invariant}\label{S:lift}
Let $X$ be a $\Z[\Z]$--homology $S^1\times S^3$, that is, a smooth
oriented $4$-manifold such that $H_*(X;\Z) = H_*(S^1\times S^3;\Z)$
and $H_*(\tilde X;\Z) = H_*(S^3;\Z)$, where $\tilde X$ is the
universal abelian cover of $X$. Choose an embedded connected
3--manifold $Y \subset X$ whose fundamental class generates
$H_3(X;\Z) = \Z$ (note that $Y$ need not be a homology sphere).
Define the Rohlin invariant of $X$ by the formula
\[
\rho(X) = \rho(Y,\sigma)\pmod 2.
\]
Here, $\sigma$ is a spin structure on $Y$ induced by either of the
two spin structures on $X$ (the two choices differ by a cohomology
class that vanishes on $Y$ hence induce the same spin structure on
$Y$), and $\rho(Y,\sigma) = \sign(W)/8\pmod 2$ is the usual Rohlin
invariant, where $W$ is any smooth compact spin $4$-manifold with
spin boundary $(Y,\sigma)$. The invariant $\rho(X)$ is well--defined,
see \cite{ruberman:ds}.

Let $f: X\to S^1$ be a smooth map such that $f_*: \pi_1 X\to \Z$
is onto. Let $Y \subset X$ be the preimage of a regular value of
$f$ then the fundamental class of $Y$ generates $H_3 (X;\Z)$ and
$\rho(X) = \rho(Y,\sigma)$. Cut $X$ open along $Y$ to obtain a spin
homology cobordism $U$ from $Y$ to itself. In particular, we have
\[
\tilde X = \ldots \cup U \cup U \cup \ldots
\]
Given a smooth compact spin $4$-manifold $W$ with boundary $\p W =
Y$, consider
\[
W_U = W \cup U \cup U \cup \ldots,
\]
which is a $4$-manifold with periodic end modeled on $\tilde X$. Any
choice of metric $g$ on $X$ naturally induces metrics on $\tilde X$
and on the end of $W_U$. The latter metric can be completed to a
metric on $W_U$ by choosing an appropriate metric on $W$. The
metrics on $\tilde X$ and $W_U$ will be called $g$ again.

By the excision principle, the Dirac operator $D (W_U,g)$ is
Fredholm if and only if the operator $D(\tilde X,g)$ is Fredholm.
According to Theorem \ref{T1}, for a generic metric $g$ on $X$,
the Dirac operator
\[
D (W_U,g): L^2_1\,(W_U,S) \to L_2\,(W_U,S)
\]
and its chiral counterparts
\[
D^{\pm} (W_U,g): L^2_1\,(W_U,S^{\pm}) \to L_2\,(W_U,S^{\mp})
\]
are Fredholm. With respect to this choice of metric, define the
metric dependent invariant
\[
w(X,g) = \ind_{\C} D^+(W_U,g) + \frac 1 8\,\sign (W).
\]
Observe that, in general, $w(X,g)$ is a rational number. It is an
integer if $\p W = Y$ is an integral homology sphere.

\begin{theorem}\label{T:rohlin}
The invariant $w(X,g)$ is well--defined (i.e. is independent of $W$)
and $w(X,g) = \rho(X)
\pmod 2$.
\end{theorem}

\begin{proof}
Let $W$ and $W'$ be two choices of a smooth compact spin $4$-manifold
with boundary $\Sigma$. Using the excision principle for the Dirac
operator, see Bunke \cite{bunke} and Charbonneau \cite{char}, and
the index theorem, we obtain
\begin{multline}\notag
\ind_{\C} D(W'_U) - \ind_{\C} D(W_U) = \ind_{\C} D(-W \cup W') \\
= -\frac 1{24}\,\int_{-W\cup W'} p_1 = -\frac 1 8\,\sign(-W\cup W')
= \frac 1 8\,\sign W - \frac 1 8\, \sign W',
\end{multline}
which proves that $w(X,g)$ is independent of the choice of $W$.

Let $Y$ and $Y'$ be two choices for cutting $X$ open.  While $Y$
and $Y'$ need not be disjoint in $X$, they become disjoint in a
finite cyclic covering $X_p \to X$ with sufficiently many sheets,
and thus become the boundary components of a smooth spin manifold
$V \subset X_p$. Because of the condition $H_*(\tilde X;\Z) =
H_*(S^3;\Z)$, the manifold $V$ must be a homology cobordism. The
independence of $w(X,g)$ of the choice of $Y$ now follows by the
argument in the previous paragraph applied to $W' = W \cup V$.

To verify the second claim observe that $\ind_{\C} D^+(W_U,g)$ is
always even because $D^+(W_U,g)$ is quaternionic linear. Therefore,
\[
w(X,g) = \ind_{\C} D^+(W_U,g) + \frac 1 8\,\sign W =
\frac 1 8\,\sign W = \rho(X) \pmod 2.
\]
\end{proof}


\section{Metrics of positive scalar curvature}\label{S:psc}
In this section we will study the $\beta$--invariant of Cappell
and Shaneson defined in \cite{cappell-shaneson:rp4} for certain
non-orientable $4$-manifolds (they called their invariant
$\alpha$ but we will call it $\beta$ to avoid confusion with the
$\alpha$--invariant of Section \ref{S:alpha1}). While the original
definition used the language of surgery theory and normal maps, we
will present it in perhaps simpler terms. After that, we will use
end periodic techniques to show that some manifolds with
non--trivial $\bcs$ invariant cannot admit metrics of positive
scalar curvature.

\subsection{The Cappell--Shaneson invariant}\label{S:cs}
Let $X'$ be a non-orientable manifold with $H_1(X';\Z) = \Z$, with
a chosen generator of $H^1(X';\Z)$. This generator comes from a map
$f: X' \to S^1$, and reduces mod 2 to $w_1 (X')$.   The orientation
double cover $X \to X'$ has a canonical orientation,  and we assume
in addition that $H_1(X;\Z) = \Z$ and $w_2 (X) =0$. Let $t:X \to X$
be the covering translation; it is of course orientation reversing.

Choose a regular value for $f$; its preimage in $X'$ is an embedded
$3$-manifold $Y$, which is transversally oriented.  Assume that $Y$
is connected, and cut $X'$ open along $Y$ to obtain a manifold $V$.
Note that both $V$ and $Y$ are orientable, although not necessarily
oriented.  Choose a lift of $Y$ to $X$, and note that this lift is
oriented because it is transversally oriented, and $X$ itself is
canonically oriented.  This choice determines a particular lift of
$V$ to $X$ as well: choose the lift of $V$ such that the normal
vector to the given
lift of $Y$ (determined, as noted before, by the map to the circle)
points out of $V$. Note that our lift of $Y$ acquires a spin structure
$\sigma$ from that of $X$ ($\sigma$ is independent of choice of spin
structure on $X$ since the two choices differ by a cohomology class
that vanishes on $Y$). Following~\cite{cappell-shaneson:rp4}, let
\[
\bcs(X',f) = \rho(Y,\sigma) - \frac 1 {16}\,\sign(V) \pmod 2,
\]
where we are implicitly using the orientation on the lift of $V$
specified above.


\begin{lemma}
The invariant $\bcs(X',f)$ is independent of the choice of $f$ but
may change sign for a different choice of the lift of $Y$ to $X$.
\end{lemma}

\begin{proof}
Let $Y_0$ and $Y_1$ be two choices. Note that we can choose homotopic
functions $f_k: X' \to S^1$, $k = 0, 1$, with $1 \in S^1$ a regular
value for both, and $Y_k = f_k^{-1}(1)$. Choose a homotopy $F: X'
\times I \to S^1$ between $f_0$ and $f_1$ such that $1 \in S^1$ is a
regular value for $F$. Then $W = F^{-1}(1)$ is a spin cobordism
between $Y_0$ and $Y_1$ which is transversally oriented. Note that
the choice of lift of $Y_0$ gives a particular lift of $W$ to $X
\times I$, and therefore a particular lift of $Y_1$ as well.

Observe that the union of $W$ and $t(W)$ separates $X \times I$ into
two pieces, one of which has oriented boundary $W \cup -V_1 \cup t(W)
\cup V_0$. Therefore,
\[
\sign(V_1) - \sign(V_0) = 2\sign(W)
\]
but also
\[
\rho(Y_1) - \rho(Y_0) = \frac 1 8\,\sign(W).
\]

\medskip\noindent
Plugging these two formulas into the definition completes the proof.
\end{proof}

The following example shows why the Rohlin invariant of $Y$ is, by
itself, not an invariant of $X'$.

\begin{example}\label{E:poincare}
Let $Y = \Sigma(2,3,5)$ be the Poincar\'e homology sphere then $Y$
bounds a smooth spin manifold $W$ with intersection form $E_8$ built
with only $2$--handles. It is easy to see (this is Kirby--calculus
folklore, cf.~\cite{ruberman:imbedding})
that the double of $W$ along $Y$ is diffeomorphic to $\#_8\,S^2
\times S^2$. Let $X'$ be the manifold obtained from the double of
$W$ by puncturing each copy of $W$ and gluing the resulting
3--sphere boundary components to get a non-orientable manifold; this
is of course the same as the connected sum of the double of $W$ with
$S^1\tilde\times S^3$.

We can compute $\bcs(X')$ by spllitting either along $Y$, or along
the $3$--sphere. Note that $\rho(Y) = 1\pmod 2$ while $\rho (S^3)
= 0 \pmod 2$. Computing with the splitting of $X'$ along $Y$, we get
that $\bcs(X',f) = \rho(Y) - \sign (W \# W)/16 = 1- (-1) = 0\pmod 2$.
On the other hand, splitting $X'$ along the $3$-sphere gives
$\bcs(X',f) = \rho(S^3) - \sign (\#_8 S^2 \times S^2)/16 = 0- (0) =
0\pmod 2$. Note that $X'$ has a metric of positive scalar curvature,
which is consistent with Theorem \ref{T:pscobs} below.
\end{example}

\subsection{An obstruction to positive scalar curvature}
Let $X'$ be a non--orien\-table $4$-manifold as described in the first
paragraph of Section \ref{S:cs}.

\begin{theorem}\label{T:pscobs}
If $X'$ admits a metric of positive scalar curvature then $\bcs(X',f) = 0$
for any choice of $f$.
\end{theorem}

\begin{proof}
Suppose that $X'$ has a metric $g$ of positive scalar curvature and define
\begin{multline}\notag
\qquad w_{\cs}(X',f,g) = \ind_{\C} D^+ (W\cup -V\cup V\cup -V\cup\ldots, g) \\
+ \sign W/8 - \sign V/16,\qquad
\end{multline}
where $W$ is any smooth compact spin manifold with boundary $Y$. That
the Dirac operator in this formula is Fredholm follows from Theorem
\ref{T:taubes} after we observe that all the operators $D^c(X,g)$ are
invertible by (\ref{E:lich}). This can also be deduced directly from
Gromov--Lawson \cite{GL}.

Observe that $w_{\cs}(X',f,g)$ reduces modulo 2 to the Cappell--Shaneson
invariant $\bcs(X',f)$. On the other hand, we have
\begin{alignat*}{1}
&w_{\cs}(X',f,g) = \\
&\ind D^+(W\cup (-V\cup V)\cup (-V\cup V)\cup\ldots)  + \sign W/8
- \sign V/16  =\\
&\ind D^+((W\cup -V)\cup (V\cup -V)\cup\ldots) + \sign(W\cup -V)/8
+ \sign V/16 = \\
&\ind D^+(-W\cup (V\cup -V)\cup \ldots) + \sign(-W)/8 + \sign V/16 =\\
&\ind D^-(W\cup (-V\cup V)\cup\ldots) + \sign(-W)/8 + \sign V/16 =\\
&- \ind D^+(W\cup (-V\cup V)\cup\ldots) - \sign W/8 + \sign V/16 =
- w_{\cs}(X',f,g)
\end{alignat*}
(the third line above is obtained by rearranging,  and the fourth by
replacing $W \cup -V$ with $-W$ and using the excision principle for
the index of the Dirac operator as in the proof of Theorem \ref{T:rohlin}).
Therefore, $w_{\cs}(X',f,g) = 0$ and $\bcs(X',f) = 0 \pmod 2$.
\end{proof}

Examples of manifolds to which Theorem~\ref{T:pscobs} applies are not
difficult to come by, and indeed are among the first known examples of
exotic $4$-manifolds. Let $Z' = S^1\tilde\times S^3$ be the
non-orientable $S^3$--bundle over $S^1$, and consider the manifold
$Z'_k = Z' \#_{k}\,(S^2 \times S^2)$. Akbulut has
constructed~\cite{akbulut:fake} a manifold $X'_1$ homotopy equivalent
(in fact, using Freedman's work, homeomorphic) to  $Z'_1$ but not
diffeomorphic to it, which is detected by the invariant $\bcs$. Using
the fact that $\beta$ is not changed by connected sum with $S^2\times
S^2$, it follows that there are such manifolds  $X'_k$  for all $k\ge
1$. According to Theorem \ref{T:pscobs} none of these manifolds admits
a metric of positive scalar curvature. The question of the existence
of an exotic $S^1\tilde\times S^3$ is still open.

Akbulut's construction is fairly strenuous, because he is trying to
keep $k$ small; there are easier constructions that produce exotic
$X'_k$ for somewhat larger $k$.

\begin{example}
The manifold $X'_{11} = Z'\#\,K3$ is homeomorphic to $Z'_{11}$. To 
see this, one can use Freedman's classification of simply connected
topological 4-manifolds to deduce that $K3$ is homeomorphic to $E\,
\#\,E\,\#_3\,S^2\times S^2$, where $E$ has the $E_8$ intersection 
form. Then 
\[
Z'\,\#\,K3 = Z'\,\#\,E\,\#E\,\#_3\,S^2\times S^2 = Z'\,\#\,(-E)\,\#
\,E\#_3\,S^2 \times S^2,
\]
where the last homeomorphism is obtained by sending $E$ around an 
orientation reversing loop. But $(-E)\,\#E\,\#_3\,S^2\times S^2$ is 
just a connected sum $\#_{11}\, S^2 \times S^2$.

The Cappell--Shaneson invariant of $X'_{11}$ equals $\rho(S^3,\sigma) 
- \sign(K3)/16 = 1\pmod 2$, therefore, $X'_{11}$ has no metric of 
positive scalar curvature.

On the other hand, the double cover of $X'_{11}$ is $(S^1\times S^3)
\,\#\,K3\,\#(-K3)$ which is diffeomorphic to $(S^1 \times S^3) \#_{22}\,
(S^2 \times S^2)$, as in Example~\ref{E:poincare}. So the double
cover of $X'_{11}$ has a metric of positive scalar curvature even though
$X'_{11}$ itself does not. Such a phenomenon has already been observed in
higher dimensions by B\'erard~Bergery~\cite{berard-bergery:isometry} and
Rosenberg~\cite{rosenberg:psc-novikov-II}, and in dimension $4$ by
LeBrun~\cite{lebrun:cover} and Hanke, Kotschick and Wehrheim \cite{HKW}. 

\end{example}


\end{document}